\documentclass{amsart}

\newtheorem{theorem}{Theorem}[section]
\newtheorem{lemma}[theorem]{Lemma}

\theoremstyle{definition}
\newtheorem{definition}[theorem]{Definition}

\newtheorem{corollary}[theorem]{Corollary}

\theoremstyle{remark}
\newtheorem{remark}[theorem]{Remark}

\numberwithin{equation}{section}

\begin{document}

\title[Stable minimal surfaces in 4-manifolds with PIC]{On complete stable minimal surfaces in 4-manifolds with positive isotropic curvature}

\author[Martin Li]{Martin Man-chun Li}
\address{Mathematics Department \\ Stanford University \\ Stanford, CA 94305}
\email{martinli@stanford.edu}


\begin{abstract}
We prove the nonexistence of stable immersed minimal surfaces uniformly conformally equivalent to $\mathbb{C}$ in any complete orientable four-dimensional  Riemannian manifold with uniformly positive isotropic curvature. We also generalize the same nonexistence result to higher dimensions provided that the ambient manifold has uniformly positive complex sectional curvature. The proof consists of two parts, assuming an ``eigenvalue condition'' on the $\overline{\partial}$-operator of a holomorphic bundle, we prove (1) a vanishing theorem for these holomorphic bundles on $\mathbb{C}$; (2) an existence theorem for holomorphic sections with controlled growth by H\"{o}rmander's weighted $L^2$-method.
\end{abstract}

\maketitle

\section{Motivation}

It has been a central theme in Riemannian geometry to study how the curvatures of a manifold affect its topology. For 2-dimensional surfaces, the Gauss-Bonnet theorem plays a key role. In particular, it implies that any oriented closed surface with positive sectional  curvature is diffeomorphic to the 2-sphere. In higher dimensions, Synge theorem says that any closed even-dimensional oriented Riemannian manifold $M^{2n}$ with positive sectional curvature is simply connected. The key idea in the proof is that there exists no stable closed geodesic $\sigma$ in such manifold. To see this, recall the second variation formula for lengths of a 1-parameter family of closed curves $\{\sigma_t\}_{t \in (-\epsilon,\epsilon)}$:
\[ \left. \frac{d^2 L(\sigma_t)}{dt^2} \right|_{t=0}=\int_{\sigma} [ \| \nabla_{\sigma'} X \|^2 - \langle R(\sigma',X)\sigma',X \rangle ] \; ds \]
where $\sigma=\sigma_0$ is a closed geodesic in $M$, $X=\left. \frac{\partial \sigma_t}{\partial t} \right|_{t=0}$ is a variation field along $\sigma$ and $R(X,Y)Z=\nabla_Y \nabla_X Z - \nabla_X \nabla_Y Z + \nabla_{[X,Y]}Z$ is the Riemann curvature tensor for $M$. Synge observed that if $M$ is even-dimensional and oriented, then there is a non-zero parallel vector field $X$ (i.e. $\nabla_{\sigma'} X=0$), hence positivity of sectional curvature implies that 
\[ \left. \frac{d^2 L(\sigma_t)}{dt^2} \right|_{t=0}<0. \]
Consequently, any closed geodesic in an oriented Riemannian manifold $M^{2n}$ with positive sectional curvature is \emph{unstable}. Therefore, $M$ has to be simply connected. Otherwise, one could minimize the length in any non-trivial free homotopy class in $\pi_1(M)$  to get a stable closed geodesic, which is a contradiction.

The idea of Synge can be extended in various directions. For example, if one considers the second variation formula for area of minimal surfaces, the curvature term that appears in the formula is the ``\emph{isotropic curvature}''. Just like positive sectional curvature tends to make geodesics unstable; positive isotropic curvature tends to make minimal surfaces unstable. In \cite{Micallef-Moore88}, M. Micallef and J. Moore proved that any minimal 2-sphere in a manifold with positive isotropic curvature is unstable. It is natural to look at the stability of complete non-compact minimal surface instead. In this paper, we prove (Theorem 2.6) the instability of complete minimal surfaces uniformly conformally equivalent to $\mathbb{C}$ by constructing holomorphic sections (half-parallel sections) with slow growth using H\"{o}rmander's weighted $L^2$-method and then applying a weighted second variation argument.

\section{Definitions and preliminaries}

The purpose of this paper is to study complete minimal surfaces $\Sigma$ in an $n$-dimensional Riemannian manifold $M$ ($n \geq 4$) which minimize area up to second order. In particular, we prove that if $n=4$ and the ambient manifold $M$ is orientable and has uniformly positive isotropic curvature, then there does not exist a complete stable minimal surface which is uniformly conformally equivalent to the complex plane $\mathbb{C}$. Note that we do not require $M$ to be compact. Therefore the result applies to universal covers of compact manifolds with positive isotropic curvature. At the end of the paper, we prove that the same result holds for any $n \geq 4$ but $M$ satisfies the stronger condition that it has uniformly positive complex sectional curvature, where $M$ need not be orientable.

First we recall some definitions. Let $M^n$ be an $n$-dimensional Riemannian manifold. 
Consider the complexified tangent bundle $T^\mathbb{C}M=TM \otimes_\mathbb{R} \mathbb{C}$ 
equipped with the Hermitian extension $\langle \cdot, \cdot \rangle$ of the inner 
product on $TM$, the curvature tensor extends to complex vectors by linearity, and 
the \emph{complex sectional curvature} of a complex two-dimensional subspace $\pi$ of $T^\mathbb{C}_pM$ at some point $p \in M$ is defined by $K^\mathbb{C}(\pi)=\langle R(v,w)\overline{w},v\rangle$, where 
$\{v,w\}$ is any unitary basis for $\pi$. 

\begin{definition}
A Riemannian manifold $M^n$ has \emph{uniformly positive complex sectional curvature} if there exists a constant $\kappa>0$ such that 
$K^\mathbb{C}(\pi) \geq \kappa >0$ for every complex two-dimensional subspace $\pi$ in $T^\mathbb{C}_pM$ at every $p \in M$.
\end{definition}

\begin{remark}
It is clear that having uniformly positive complex sectional curvature implies having uniformly positive sectional curvature. (One simply considers all $\pi$ which comes from the complexification of a real two-dimensional subspace in $T_pM$.) Therefore, by Bonnet-Myers theorem, $M$ is automatically compact if it is complete.
\end{remark}

Using Ricci flow techniques, S. Brendle and R. Schoen \cite{Brendle-Schoen09} proved that any manifold $M$ with uniformly positive complex sectional curvature is diffeomorphic to a spherical space form. In particular, when $M$ is simply-connected, $M$ is diffeomorphic to the $n$-dimensional sphere. In fact, Brendle and Schoen showed that the condition of positive complex sectional curvature is preserved under the Ricci flow, and any manifold equipped with such a metric evolves under the normalized Ricci flow to a spherical space-form. The optimal convergence result so far is obtained by S. Brendle in \cite{Brendle08}, where he proved that any compact manifold $M$ such that $M \times \mathbb{R}$ has positive isotropic curvature would converge to a spherical space form under the normalized Ricci flow.

As a result, having positive complex sectional curvature is a rather restrictive condition. There is a related positivity condition which allows more flexibility called \emph{positive isotropic curvature} (PIC). Instead of extending in a Hermitian way, one can choose to extend the metric on $TM$ to a $\mathbb{C}$-bilinear form $(\cdot,\cdot)$ on $T^\mathbb{C}M$. We say that a vector $v \in T^\mathbb{C}_pM$ is \emph{isotropic} if $(v,v)=0$. A subspace $V \subset T^\mathbb{C}_p$ is \emph{isotropic} if every $v \in V$ is isotropic. 

\begin{definition}
A Riemannian manifold $M^n$, $n \geq 4$, has \emph{uniformly positive isotropic curvature} (uniformly PIC) if there exists a positive constant $\kappa>0$ such that $K^\mathbb{C}(\pi) \geq \kappa >0$ for every isotropic complex two-dimensional subspace $\pi \subset T^\mathbb{C}_pM$ at every $p \in M$.
\end{definition}

\begin{remark}
This condition is clearly weaker than the previous one because it only requires positivity of complex sectional curvature among \emph{isotropic} 2-planes. Moreover, the condition is vacuous when $n\leq 3$ (see \cite{Micallef-Moore88}), so we restrict ourselves to $n \geq 4$.
\end{remark}

In \cite{Brendle-Schoen09}, it was shown that uniformly PIC is also a condition preserved under the Ricci flow. However, it is not true that any manifold equipped with such a metric would converge to a spherical space-form under the normalized Ricci flow. For example, the product manifold $S^{n-1} \times S^1$ is PIC but it has no metric with constant positive sectional curvature. (The universal cover of $S^{n-1} \times S^1$ is $S^{n-1} \times \mathbb{R}$, which is not $S^n$.)

On the other hand, there is a ``sphere theorem'' for compact manifolds with PIC. Namely, if $M^n$, $n\geq 4$, is a compact manifold with PIC, and $M$ is simply-connected, then $M^n$ is \emph{homeomorphic} to $S^n$. This sphere theorem, which generalizes the classical sphere theorem since any $1/4$-pinched manifold is PIC, was proved by M. Micallef and J. Moore in \cite{Micallef-Moore88}, where the notion of PIC was first defined. In other words, all simply-connected PIC manifolds are topologically ``trivial''. This, of course, raises the natural question: what are the topological obstructions on the fundamental group $\pi_1M$ for a compact manifold $M$ to admit a metric with PIC?

Along this direction, A. Fraser \cite{Fraser03} proved that for $n \geq 5$, $\pi_1M$ cannot contain a subgroup isomorphic to $\mathbb{Z} \oplus \mathbb{Z}$, the fundamental group of a torus. Her proof uses the existence theory of stable minimal surfaces by R. Schoen and S.T. Yau \cite{Schoen-Yau79}, and the Riemann-Roch theorem to construct \emph{almost holomorphic} variations which, in turn, contradicts stability. A few years later, S. Brendle and R. Schoen \cite{Brendle-Schoen} proved that the same result holds for $n=4$. On the other hand, M. Micallef and M. Wang \cite{Micallef-Wang93} proved that if $(M_1^n,g_1)$ and $(M_2^n,g_2)$ are manifolds with positive isotropic curvature, then $M_1 \sharp M_2$ also admits a metric with positive isotropic curvature. In particular, $\sharp_{i=1}^k(S^1 \times S^{n-1})$ admits a metric with positive isotropic curvature for any positive integer $k$. Therefore, the fundamental group of a manifold with positive isotropic curvature can be quite large.

To state the main result in this paper, we need the following definition.

\begin{definition}
A Riemannian surface $(\Sigma,h)$ is \emph{uniformly conformally equivalent} to the complex plane $\mathbb{C}$ if 
there is a (conformal) diffeomorphism $\phi:\mathbb{C} \to \Sigma$, a positive smooth function $\lambda$ on $\mathbb{C}$ 
and a constant $C>0$ such that 
\begin{equation}
\phi^*h=\lambda^2 |dz|^2 \qquad \text{with} \qquad \frac{1}{C} \leq \lambda^2. 
\end{equation}
\end{definition}

Our theorem states that surfaces of this type cannot arise as stable minimal surfaces in any orientable 4-manifold with uniformly positive isotropic curvature.

\begin{theorem}
Let $M$ be a 4-dimensional complete orientable Riemannian manifold with uniformly positive isotropic curvature. Let $\mathbb{C}$ be the complex plane equipped with the standard flat metric. Then there does not exist a stable immersed minimal surface $\Sigma$ in $M$ which is uniformly conformally equivalent to $\mathbb{C}$.
\end{theorem}

If we assume that $M$ has uniformly positive complex sectional curvature, then the result holds in any dimension, and without the orientability assumption on $M$.

\begin{theorem}
Let $M$ be an $n$-dimensional complete (not necessarily orientable) Riemannian manifold with uniformly positive complex sectional curvature. Let $\mathbb{C}$ be the complex plane equipped with the standard flat metric. Then there does not exist a stable immersed minimal surface $\Sigma$ in $M$ which is uniformly conformally equivalent to $\mathbb{C}$.
\end{theorem}

An outline of the paper is as follows. In section 3, we prove a vanishing theorem for the $\overline{\partial}$-operator on holomorphic bundles satisfying an ``eigenvalue condition''. In section 4, we describe H\"{o}rmander's weighted $L^2$-method and use it to construct holomorphic sections with controlled growth, assuming the ``eigenvalue condition''. In section 5, we apply the results to prove Theorem 2.6 and Theorem 2.7.

\textit{Acknowledgement.} The author would like to thank his advisor, Professor Richard Schoen, for all his continuous support and encouragement throughout the progress of this work. He would also like to thank Professor Simon Brendle for many helpful discussions. He is thankful to the referee for many useful comments.

\section{A vanishing theorem}
 
Throughout the paper, $\mathbb{C}$ will denote the standard complex plane with the flat metric 
\begin{equation*} 
ds^2=dx^2+dy^2=|dz|^2. 
\end{equation*}
Suppose $E$ is a holomorphic vector bundle over $\mathbb{C}$ with a compatible Hermitian metric. Let $\overline{\partial}$ denote the $\overline{\partial}$-operator associated to $E$. Assume that $(E,\overline{\partial})$ satisfies the following ``eigenvalue condition'': these exists some positive constant $\kappa_0>0$, and a sufficiently small constant $\epsilon_0>0$ (depending on $\kappa_0$) such that for all $0< \epsilon<\epsilon_0$, we have  
\begin{equation}
\kappa_0 \int_\mathbb{C} |s|^2 e^{-\epsilon |z|} \; dxdy \leq \int_\mathbb{C} |\overline{\partial}  s|^2 e^{-\epsilon |z|} \; dxdy
\end{equation}
for all compactly supported smooth sections $s$ of $E$, we write $s \in C^\infty_c(E)$.

For any positive continuous function $\varphi$ on $\mathbb{C}$, we can define an $L^2$-norm on $C^\infty_c(E)$ by
\begin{equation*}
\|s\|_{L^2(E,\varphi)}=\left( \int_\mathbb{C} |s|^2 \varphi \; dxdy \right)^{\frac{1}{2}} .
\end{equation*}
Let $L^2(E,\varphi)$ be the Hilbert space completion of $C^\infty_c(E)$ with respect to the weighted $L^2$-norm $\| \cdot \|_{L^2(E,\varphi)}$. In other words, $s \in L^2(E,\varphi)$  if and only if $s$ is a measurable section of $E$ with $\|s\|_{L^2(E,\varphi)}<\infty$.

In this section, we prove a vanishing theorem for the $\overline{\partial}$-operator on these weighted $L^2$ spaces of sections of holomorphic bundles over $\mathbb{C}$ satisfying the ``eigenvalue condition'' above. Roughly speaking, a holomorphic section cannot grow too slowly unless it is identically zero.

\begin{theorem}
Suppose $E$ is a holomorphic vector bundle over $\mathbb{C}$ satisfying the ``eigenvalue condition'' (3.1) with constants $\kappa_0$ and $\epsilon_0$. Then, there exists no non-trivial holomorphic section $s \in L^2(E,e^{-4\epsilon |z|})$ for any $0 < \epsilon<\epsilon_0/4$.
\end{theorem}

\begin{proof}
The proof is a direct cutoff argument. Suppose $s$ is a holomorphic section of $E$ which belongs to $L^2(E,e^{-4\epsilon |z|})$ for some $0< \epsilon<\epsilon_0/4$. We will show that $s \equiv 0$.

For every real number $R>0$, choose a cutoff function $\phi_R \in C^\infty_c(\mathbb{C})$ such that 
\begin{itemize}
	\item $\phi_R(z)=1$  for $|z| \leq R$;
	\item $\phi_R(z)=0$ for $|z| \geq 2R$;
	\item $|\nabla \phi_R| \leq \frac{2}{R}$.
\end{itemize}
Let $\hat{s}=\phi_R s$. Note that $\hat{s} \in C^\infty_c(E)$. Therefore, by (3.1), properties of $\phi_R$ and holomorphicity of $s$,
\begin{align*} 
\kappa_0 \int_{|z| \leq R} |s|^2 e^{-4\epsilon|z|} \; dxdy & \leq \kappa_0 \int_\mathbb{C} |\hat{s}|^2 e^{-4\epsilon|z|} \; dxdy \\
&\leq \int_\mathbb{C} |\overline{\partial}\hat{s}|^2 e^{-4\epsilon |z|} \; dxdy \\
&=\int_\mathbb{C} |\overline{\partial}\phi_R|^2 |s|^2 e^{-4\epsilon|z|} \; dxdy \\
& \leq \frac{1}{R^2} \int_{R \leq |z| \leq 2R} |s|^2 e^{-4\epsilon |z|} \; dxdy  \\
&\leq \frac{1}{R^2} \|s\|^2_ {L^2(E,e^{-4\epsilon |z|})}.
\end{align*}
By our assumption, $s \in L^2(E,e^{-4\epsilon |z|})$. As $R \to \infty$, the right hand side goes to zero while the left hand side goes to $\kappa_0 \|s\|_{L^2(E,e^{-4\epsilon |z|})}$. Since $\kappa_0>0$, we conclude that $\|s\|_{L^2(E,e^{-4\epsilon |z|})}=0$, hence $s \equiv 0$. This completes the proof.
\end{proof}

\section{H\"{o}rmander's weighted $L^2$ method}

In this section, we will use H\"{o}rmander's weighted $L^2$ method to construct non-trivial weighted $L^2$ holomorphic sections on $E$. Some basic facts on unbounded operators between Hilbert spaces can be found in \cite{Chen-Shaw}. 

Assume $E$ is a holomorphic vector bundle on $\mathbb{C}$ satisfying the ``eigenvalue condition'' (3.1) with constants $\epsilon_0$ and $\kappa_0$. For this section, we also assume that $E$ is the complexification of some real vector bundle $\xi$, hence $\overline{E}=E$.

Recall that there is a natural first order differential operator $\overline{\partial}$ defined on the space of smooth compactly supported sections of $E$:
\begin{equation*}
\overline{\partial}: C^\infty_c(E) \to C^\infty_c(E \otimes T^{0,1}\mathbb{C}),
\end{equation*}
where
\begin{equation*}
\overline{\partial}s=(\nabla_{\frac{\partial}{\partial \overline{z}}} s )\otimes d\overline{z}.
\end{equation*}
Let $L^2(E,e^{-2\epsilon |z|})$ be the Hilbert space completion of $C^\infty_c(E)$ with respect to the weighted $L^2$-norm
\begin{equation*} 
\|s\|_{L^2(E,e^{-2\epsilon |z|})} = \left( \int_\mathbb{C} |s|^2 e^{-2 \epsilon |z|} \; dxdy \right)^{\frac{1}{2}}. 
\end{equation*}
Similarly, let $L^2(E \otimes T^{0,1} \mathbb{C}, e^{-2\epsilon |z|})$ denote the Hilbert space completion of $C^\infty_c(E \otimes T^{0,1}\mathbb{C})$ with respect to the weighted $L^2$-norm 
\begin{equation*} 
\| \sigma \|_{L^2(E \otimes T^{0,1}\mathbb{C},e^{-2\epsilon |z|})} = \left( \int_\mathbb{C} |\sigma|^2 e^{-2 \epsilon |z|} \; dxdy \right)^{\frac{1}{2}} . 
\end{equation*}
Now, let 
\begin{equation*} 
\overline{\partial}: L^2(E,e^{-2\epsilon |z|}) \to L^2(E \otimes T^{0,1}\mathbb{C}, e^{-2 \epsilon |z|})
\end{equation*}
be the maximal closure of $\overline{\partial}$ defined as follows: an element $s \in L^2(E,e^{-2\epsilon |z|})$ is in the domain of $\overline{\partial}$ if $\overline{\partial}s$, defined in the distributional sense, belongs to $L^2(E \otimes T^{0,1}\mathbb{C},e^{-2\epsilon |z|})$. Then, $\overline{\partial}$ defines a linear, closed, densely defined unbounded operator. Note that $\overline{\partial}$ is closed because differentiation is a continuous operation in distribution theory. It is densely defined since Dom($\overline{\partial}$) contains all compactly supported smooth sections $C^\infty_c(E)$, which is clearly dense in $L^2(E,e^{-2\epsilon |z|})$.

By standard Hilbert space theory, the Hilbert space adjoint of $\overline{\partial}$, denoted by $\overline{\partial}^*$, is a linear, closed, densely defined unbounded operator and 
\begin{equation*} 
\overline{\partial}^* : L^2(E \otimes T^{0,1}\mathbb{C},e^{-2\epsilon |z|}) \to L^2(E,e^{-2 \epsilon |z|}). 
\end{equation*}
An element $\sigma$ belongs to Dom$(\overline{\partial}^*)$ if there is an $s \in L^2(E,e^{-2 \epsilon |z|})$ such that for every $t \in $Dom$(\overline{\partial})$, we have 
\begin{equation*} 
(\sigma,\overline{\partial}t)_{L^2(E \otimes T^{0,1}\mathbb{C},e^{-2\epsilon |z|})}=(s,t)_{L^2(E,e^{-2\epsilon |z|})}. 
\end{equation*}
We then define $\overline{\partial}^* \sigma=s$. 

The key theorem in this section is the surjectivity of $\overline{\partial}$.

\begin{theorem}
Suppose $E$ is a holomorphic vector bundle over $\mathbb{C}$ satisfying the ``eigenvalue condition'' (3.1) with constants $\kappa_0$ and $\epsilon_0$. Moreover, assume that $E$ is the complexification of some real vector bundle $\xi$ over $\mathbb{C}$. Then,
\begin{equation*} 
\overline{\partial} : L^2(E,e^{-2\epsilon |z|}) \to L^2(E \otimes T^{0,1}\mathbb{C}, e^{-2 \epsilon |z|}) 
\end{equation*}
is surjective for all $0 < \epsilon < \min (\frac{\epsilon_0}{2},\frac{\sqrt{\kappa_0}}{2})$.
\end{theorem}

\begin{proof}
By standard Hilbert space theory (see section 4.1 in \cite{Chen-Shaw}), it suffices to show that 
\begin{itemize}
	\item[(i)] the adjoint operator $\overline{\partial}^*$ is injective, and 
	\item[(ii)] the range of $\overline{\partial}$ is closed.
\end{itemize}

First of all, we need to compute $\overline{\partial}^*$ explicitly.

\textit{Claim 1:} For any $\sigma=s \otimes d\overline{z} \in C^\infty_c(E \otimes T^{0,1}\mathbb{C})$,
\begin{equation} 
\overline{\partial}^* \sigma = - \nabla_{\frac{\partial}{\partial z}} s + \epsilon \frac{\overline{z}}{|z|} s.
\end{equation}

\textit{Proof of Claim 1:} This is just integration by parts. Let $t \in $Dom$(\overline{\partial})$, and $\langle \cdot, \cdot \rangle_E$, $\langle \cdot, \cdot \rangle_{E \otimes T^{0,1}\mathbb{C}}$  be the pointwise Hermitian metric on $E$ and $E \otimes T^{0,1}\mathbb{C}$ respectively.
\begin{align*}
(\sigma,\overline{\partial} t)_{L^2(E \otimes T^{0,1}\mathbb{C},e^{-2\epsilon |z|})}&=\int_\mathbb{C} \langle \sigma, \overline{\partial} t \rangle_{E \otimes T^{0,1}\mathbb{C}} \;e^{-2 \epsilon |z|}\; dxdy \\
&= \int_\mathbb{C} \langle  e^{-2 \epsilon |z|} \sigma, \overline{\partial} t \rangle_{E \otimes T^{0,1}\mathbb{C}}\; dxdy 
\\
&= \int_\mathbb{C} \langle  e^{-2 \epsilon |z|} s, \nabla_{\frac{\partial}{\partial \overline{z}}} t \rangle_E \; dxdy 
\\
&= \int_\mathbb{C} \frac{\partial}{\partial z} \langle e^{-2\epsilon |z|} s,t \rangle_E \; dxdy -
\int_\mathbb{C} \langle \nabla_{\frac{\partial}{\partial z}}(e^{-2 \epsilon |z|} 
s),t \rangle_E \; dxdy \\
&= -\int_\mathbb{C} \langle \nabla_{\frac{\partial}{\partial z}} s -\epsilon \frac{\overline{z}}{|z|}s,t 
\rangle_E \; e^{-2\epsilon |z|} \; dxdy. \\
&=(-\nabla_{\frac{\partial}{\partial z}} s + \epsilon \frac{\overline{z}}{|z|} s,t)_{L^2(E,e^{-2\epsilon |z|})}.
\end{align*}
The first term in the second to last line vanishes because we have integrated by parts and used the fact that $s$ is compactly supported. This proves Claim 1.

Next, we need to establish a basic estimate for the adjoint operator $\overline{\partial}^*$.

\textit{Claim 2:} For every $0 < \epsilon < \min(\frac{\epsilon_0}{2},\frac{\sqrt{\kappa_0}}{2})$, there exists a constant $\kappa_1>0$ such that 
\begin{equation}  
\kappa_1 \int_\mathbb{C} |\sigma|^2 e^{-2 \epsilon |z|} \; dxdy \leq \int_\mathbb{C} | \overline{\partial}^* \sigma|^2 e^{-2 \epsilon |z|} \; dxdy 
\end{equation}
for all $\sigma \in $Dom($\overline{\partial}^*)$.

\textit{Proof of Claim 2:} Since $E$ is the complexification of some real vector bundle, (2.1) implies that
\begin{equation} 
\kappa_0 \int_\mathbb{C} |s|^2 e^{-2\epsilon |z|} \; dxdy \leq \int_\mathbb{C} |\nabla_{\frac{\partial}{\partial z}} s |^2 e^{-2\epsilon |z|} \; dxdy 
\end{equation}
for any $s \in C^\infty_c(E)$. First, we establish claim 2 for $\sigma=s \otimes d\overline{z} \in C^\infty_c(E \otimes T^{0,1}\mathbb{C})$. By (3.1), the triangle inequality, and that $\epsilon < \frac{\sqrt{\kappa_0}}{2}$, we have
\begin{align*}
\| \overline{\partial}^* \sigma\|_{L^2(E,e^{-2\epsilon |z|})}
&=\| -\nabla_{\frac{\partial}{\partial z}} s + \epsilon \frac{\overline{z}}{|z|} s \|_{L^2(E,e^{-2\epsilon |z|})} \\
&\geq \| \nabla_{\frac{\partial}{\partial z}} s\|_{L^2(E,e^{-2\epsilon |z|})} - \epsilon \| s \|_{L^2(E,e^{-2\epsilon |z|})} \\
&\geq \sqrt{\kappa_0} \|s \|_{L^2(E,e^{-2\epsilon |z|})} - \epsilon \| s \|_{L^2(E,e^{-2\epsilon |z|})} \\
&\geq \frac{\sqrt{\kappa_0}}{2} \|s \|_{L^2(E,e^{-2\epsilon |z|})} \\
&=  \frac{\sqrt{\kappa_0}}{2} \|\sigma \|_{L^2(E \otimes T^{0,1}\mathbb{C},e^{-2\epsilon |z|})}. 
\end{align*} 
Squaring both sides give the inequality we want, with $\kappa_1=\frac{\kappa_0}{4}$. To prove the inequality for arbitrary $\sigma \in $Dom($\overline{\partial}^*)$, it suffices to show that $C^\infty_c(E\otimes T^{0,1}\mathbb{C})$ is dense in Dom($\overline{\partial}^*$) in the graph norm
\begin{equation*} 
\sigma \mapsto \|\sigma\|_{L^2(E \otimes T^{0,1}\mathbb{C},e^{-2\epsilon |z|})} + \| \overline{\partial}^* \sigma \|_{L^2(E,e^{-2\epsilon |z|})}. 
\end{equation*}
A proof of this elementary fact can be found in Appendix A. Therefore, we have completed the proof of claim 2.

Now, claim 2 clearly implies both (i) and (ii) (lemma 4.1.1 in \cite{Chen-Shaw}). This finishes the proof of Theorem 4.1.
\end{proof}

An important corollary of Theorem 4.1 is the following existence theorem of holomorphic sections of $E$ with controlled growth.

\begin{corollary}
Suppose $E$ is a holomorphic vector bundle over $\mathbb{C}$ satisfying the ``eigenvalue condition'' (3.1) with constants $\kappa_0$ and $\epsilon_0$. Assume that $E$ is the complexification of some real vector bundle $\xi$ over $\mathbb{C}$.  

Then, for any $0 < \epsilon <\min (\frac{\epsilon_0}{2},\frac{\sqrt{\kappa_0}}{2})$, there exists a non-trivial holomorphic section $s \in L^2(E,e^{-4\epsilon |z|})$, that is, 
\begin{equation*} 
\overline{\partial}s=0 
\end{equation*}
and 
\begin{equation*} 
\int_\mathbb{C} |s|^2 e^{-4 \epsilon |z|} \; dxdy < \infty. 
\end{equation*}
\end{corollary}

\begin{proof}
First, notice that any holomorphic vector bundle on a non-compact Riemann surface is 
holomorphically trivial (\cite{Forster}). Therefore, we can choose a nowhere vanishing holomorphic 
section $u$ of $E$. However, such a $u$ maybe not be in $L^2(E,e^{-4\epsilon |z|})$. We will correct it by a 
cutoff argument and solving an inhomogeneous equation of the form 
$\overline{\partial}s=\sigma$ to construct a holomorphic section in 
$L^2(E,e^{-4\epsilon |z|})$. 

Take a smooth compactly supported cutoff function $\psi \in C^\infty_c(\mathbb{C})$ such that 
\begin{itemize}
\item $\psi(z)=1$ for $|z| \leq 1$, and 
\item $\psi(z)=0$ for $|z| \geq 2$.
\end{itemize}

Define 
\begin{equation*} 
v=\frac{1}{z} \left(\frac{\partial \psi}{\partial \overline{z}} \right) u. 
\end{equation*}
Observe that $v \in C^\infty_c(E)$ even though $1/z$ is singular at $z=0$. 
By Theorem 4.1, there exists $w \in L^2(E,e^{-2\epsilon |z|})$ such that 
\begin{equation*} 
\overline{\partial}w=v \otimes d\overline{z} .
\end{equation*}
Since $v$ is smooth and compactly supported, by elliptic regularity, $w$ is a smooth section of $E$ (but not necessarily compactly supported). 

Next, we let
\begin{equation*} 
s=\psi u-zw. 
\end{equation*}

\textit{Claim:} $s$ is a non-trivial holomorphic section in $L^2(E,e^{-4\epsilon |z|})$. 

\textit{Proof of Claim:} First of all, 
\begin{equation*} 
s(0)=\psi(0)u(0)=u(0) \neq 0 
\end{equation*}
since $u$ is nowhere vanishing. $s$ is holomorphic since 
\begin{equation*} 
\overline{\partial}s=\overline{\partial} (\psi u)-z\overline{\partial}w
=\frac{\partial \psi}{\partial \overline{z}} u \otimes d\overline{z} - z \overline{\partial}w=0. 
\end{equation*}
Finally, to see that $s$ is in $L^2(E,e^{-4\epsilon |z|})$, we see that $\psi u$ is smooth and compactly supported, 
hence, in $L^2(E,e^{-4\epsilon |z|})$. Moreover, 
\begin{equation*} 
\int_\mathbb{C} |z|^2 |w|^2 e^{-4\epsilon |z|} \; dxdy = \int_\mathbb{C} |z|^2e^{-2\epsilon|z|} 
|w|^2 e^{-2\epsilon |z|} \; dxdy 
\end{equation*}
Since $|z|^2e^{-2\epsilon |z|} \leq 1$ on $|z| \geq R$ for some $R>0$ sufficiently large, it follows that 
$zw \in L^2(E,e^{-4\epsilon |z|})$, and so does $s$. This proves our claim and hence establishes the corollary.
\end{proof}

\section{Proof of Theorem 2.6 and 2.7}

In this section, we prove Theorems 2.6 and 2.7 using the results in section 3 and 4. 

\textit{Proof of Theorem 2.6:} 
We argue by contradiction. Suppose Theorem 2.6 is false. Then there exists a stable minimal immersion $u: \Sigma \to M$ into an oriented Riemannian 4-manifold $M$ with uniformly positive isotropic curvature bounded from below by $\kappa>0$, with $\Sigma$ uniformly conformally equivalent to $\mathbb{C}$. Recall that $u$ is \emph{minimal} if it is a critical point of the area functional with respect to compactly supported variations, and $u$ is \emph{stable} if and only if the second variation of area for any compactly supported variation is nonnegative. Note that we assume $u$ to be an immersion, therefore we do not allow the existence of branch points.

We consider the normal bundle of the surface $u(\Sigma)$. We denote by $F$ the bundle on $\Sigma$ given by the pullback under $u$ of the normal bundle of $u(\Sigma)$. Note that $F$ is a smooth real vector bundle of rank 2. Since $M$ and $\Sigma$ are orientable, we conclude that $F$ is orientable. Let $F^\mathbb{C}=F \otimes_\mathbb{R} \mathbb{C}$ be the complexification of $F$. Since $F$ is orientable, the complexified bundle $F^\mathbb{C}$ splits as a direct sum of two holomorphic line bundles $F^{1,0}$ and $F^{0,1}$. Here, $F^{1,0}$ consists of all vectors of the form $\mu(v-iw) \in F^\mathbb{C}$, where $\mu \in \mathbb{C}$ and $\{v,w\}$ is a positively oriented orthonormal basis of $F$. An important observation here is that every section $s \in C^\infty(F^{1,0})$ is automatically isotropic, i.e. $(s,s)=0$. Moreover, the splitting $F=F^{1,0}\oplus F^{0,1}$ is parallel, i.e. invariant under $\nabla$.

Since $u:\Sigma \to M$ is stable, the complexified stability inequality (see \cite{Fraser03}) says that
\begin{equation*} 
\int_\Sigma \left\langle R\left(s,\frac{\partial u}{\partial z}\right)\frac{\partial u}{\partial \overline{z}}, s \right\rangle 
\; dxdy \leq \int_\Sigma (|\nabla^\perp_{\frac{\partial}{\partial \overline{z}}} s|^2-
|\nabla^\top_{\frac{\partial}{\partial z}} s|^2) \; dxdy 
\end{equation*}
for all compactly supported sections $s \in C^\infty_c(F)$, where $z=x+iy$ is a local isothermal coordinate on $\Sigma$. In particular, it holds for all $s \in C^\infty_c(F^{0,1})$. Note that $s \perp \frac{\partial u}{\partial z}$, $s$ is isotropic and $\frac{\partial u}{\partial z}$ is also isotropic (since $z$ is an isothermal coordinate and $u$ is an isometric immersion). Therefore, $\{s,\frac{\partial u}{\partial z}\}$ span a two dimensional isotropic subspace. Using the lower bound on the isotropic curvature and throwing away the second term on the right, we get the following inequality
\begin{equation}
\kappa \int_\Sigma |s|^2 \; da \leq \int_\Sigma |\overline{\partial}s|^2 \; da
\end{equation}
for every $s \in C^\infty_c(F^{1,0})$, where $da$ is the area element of $\Sigma$.

Since $\Sigma$ is uniformly conformally equivalent to $\mathbb{C}$, by definition, there exists a diffeomorphism $\phi:\mathbb{C} \to \Sigma$ and a constant $C>0$ such that
\begin{equation} \phi^*h=\lambda^2 |dz|^2 \qquad \text{with} \qquad \frac{1}{C} \leq \lambda^2, \end{equation}
where $h$ is the induced metric on $\Sigma$. Define $E=\phi^*(F^{1,0})$ be the pullback of the holomorphic bundle $F^{1,0}$ by $\phi$. Since $\phi$ is a conformal diffeomorphism, $E$ is again a holomorphic bundle over $\mathbb{C}$. By (5.1) and (5.2), 
\begin{equation} 
\frac{\kappa}{C} \int_\mathbb{C} |s|^2 \; dxdy \leq \int_\mathbb{C} |\overline{\partial}s|^2 \; dxdy 
\end{equation}
for every $s \in C^\infty_c(E)$. We then show that $E$ satisfies the ``eigenvalue condition'' (3.1).

\textit{Claim:} there exists a constant $\epsilon_0>0$ such that for all $0<\epsilon<\epsilon_0$, 
\begin{equation} 
\frac{\kappa}{4C} \int_\mathbb{C} |s|^2 e^{-\epsilon |z|} \; dxdy \leq \int_\mathbb{C} |\overline{\partial}  s|^2 e^{-\epsilon |z|} \; dxdy 
\end{equation}
for all compactly supported smooth sections $s \in C^\infty_c(E)$.

\textit{Proof of Claim:} Let $s \in C^\infty_c(E)$. Take $t=e^{-\epsilon |z|/2}s$, notice that 
\begin{align*}
|\overline{\partial}t|^2&=|\overline{\partial}(e^{-\epsilon |z|/2})s+e^{-\epsilon |z|/2}\overline{\partial}s|^2 \\
&\leq 2|\overline{\partial}(e^{-\epsilon |z|/2})s|^2+2|e^{-\epsilon |z|/2}\overline{\partial}s|^2 \\
&= 2 e^{-\epsilon |z|} ( \frac{\epsilon^2}{16}  |s|^2+|\overline{\partial}s|^2)
\end{align*}
Applying (5.3) to $t \in C^\infty_c(E)$ and using the above estimate, 
\begin{align*}
\frac{\kappa}{C} \int_\mathbb{C} |s|^2 e^{-\epsilon |z|} \; dxdy &= \frac{\kappa}{C} \int_\mathbb{C} |t|^2 \; dxdy \\
&\leq \int_\mathbb{C} |\overline{\partial}t|^2 \; dxdy \\
&\leq \frac{\epsilon^2}{8} \int_\mathbb{C} |s|^2 e^{-\epsilon |z|} \; dxdy+ 2 \int_\mathbb{C} |\overline{\partial}s|^2 e^{- \epsilon |z|} \; dxdy. 
\end{align*}
Hence, if we take $\epsilon_0=\frac{2 \sqrt{\kappa}}{\sqrt{C}}$, for every $0< \epsilon <\epsilon_0$, we get 
\begin{equation*} 
\frac{\kappa}{4C} \int_\mathbb{C} |s|^2 e^{-\epsilon |z|} \; dxdy \leq \int_\mathbb{C} |\overline{\partial}s|^2 e^{- \epsilon |z|} \; dxdy. 
\end{equation*}
This proves our claim.

To summarize, we have constructed a holomorphic vector bundle $E$ over $\mathbb{C}$ which satisfies the ``eigenvalue condition'' (3.1) with $\kappa_0=\frac{\kappa}{4C}$ and $\epsilon_0=\frac{2 \sqrt{\kappa}}{\sqrt{C}}$. So we can apply our results in section 3. Moreover, even though $E$ is not a complexification of a real vector bundle, we see that the result in section 4 still holds for $E$ because $E$ satisfies (4.3) (This follows from the fact that the inequality (5.1) holds with $\overline{\partial}$ replaced by $\partial$, we have to use the fact that $F=F^{1,0}\oplus F^{0,1}$ is a parallel splitting). Hence, if we fix $\epsilon>0$ sufficiently small ($\epsilon< \min(\frac{\epsilon_0}{4},\frac{\sqrt{\kappa_0}}{2})$), then Corollary 4.2 gives a non-trivial holomorphic section $s \in L^2(E,e^{-4\epsilon |z|})$, which contradicts Theorem 3.1. This contradiction completes the proof of Theorem 2.6.

\qed

\textit{Proof of Theorem 2.7:} The proof is very similar to the above. Again we argue by contradiction. Suppose Theorem 2.7 is false. Then there exists a stable minimal immersion $u: \Sigma \to M$ into a Riemannian n-manifold $M$ with uniformly positive complex sectional curvature bounded from below by $\kappa>0$, with $\Sigma$ uniformly conformally equivalent to $\mathbb{C}$. 

Let $F$ be the complexified normal bundle of $u(\Sigma)$ as before and let $E=\phi^*F$. By our assumption on $M$, (5.1) holds for every $s \in C^\infty_c(F)$. Exactly the same argument as above gives our desired contradiction.
\qed

\section{Appendix A: A density lemma}

We give a proof of the following density lemma used in the proof of Theorem 4.1. The proof is very similar to that of Lemma 4.1.3 in \cite{Hormander}.

\begin{lemma}
The subspace $C^\infty_c(E\otimes T^{0,1}\mathbb{C})$ is dense in Dom($\overline{\partial}^*$) in the graph norm
\begin{equation*} 
\sigma \mapsto \|\sigma\|_{L^2(E \otimes T^{0,1}\mathbb{C},e^{-2\epsilon |z|})} + \| \overline{\partial}^* \sigma \|_{L^2(E,e^{-2\epsilon |z|})}. 
\end{equation*}
\end{lemma}

\begin{proof}
Let $\sigma=s\otimes d\overline{z}\in $Dom$(\overline{\partial}^*)$. First of all, we will show that the set of $\tau \in $Dom$(\overline{\partial}^*)$ with compact support is dense in Dom$(\overline{\partial}^*)$. For each $R>0$, let $\varphi=\varphi_R \in C^\infty_c(\mathbb{C})$ be a smooth cutoff function such that 
\begin{itemize}
	\item $\varphi(z)=1$  for $|z| \leq R$;
	\item $\varphi(z)=0$ for $|z| \geq 2R$;
	\item $|\nabla \varphi| \leq \frac{2}{R}$.
\end{itemize}

\textit{Claim 1:} When $R \to \infty$, $\varphi_R \; \sigma$ converges to $\sigma$ in the graph norm.

\textit{Proof of Claim 1:} First of all, we observe that $\varphi_R \;  \sigma \in$Dom$(\overline{\partial}^*)$. In fact, for any $t \in $Dom$(\overline{\partial})$, 
\begin{align*}
(\varphi_R \sigma, \overline{\partial}t)
&= (\sigma, \overline{\partial}(\varphi_R t))-(\sigma,(\overline{\partial}\varphi_R)t) \\
&= (\overline{\partial}^* \sigma, \varphi_R t) - ((\frac{\partial \varphi_R}{\partial z})s,t) \\
&= (\varphi_R \overline{\partial}^* \sigma-(\frac{\partial \varphi_R}{\partial z})s,t)
\end{align*}
It follows that $(\varphi_R \sigma, \overline{\partial}t)$ is continuous in $t$ for the norm $\|t\|_{L^2(E,e^{-2\epsilon |z|})}$, so $\varphi_R \;  \sigma \in$Dom$(\overline{\partial}^*)$ and 
\begin{equation*} 
\overline{\partial}^*(\varphi_R \sigma)= \varphi_R \overline{\partial}^* \sigma-(\frac{\partial \varphi_R}{\partial z})s.
\end{equation*}
Therefore, we have the estimate
\begin{equation*} 
|\overline{\partial}^*(\varphi_R \sigma)- \varphi_R \overline{\partial}^* \sigma| \leq \frac{2}{R}|\sigma|. 
\end{equation*}
Since $\sigma \in L^2(E \otimes T^{0,1}\mathbb{C},e^{-2\epsilon |z|})$, and $\varphi_R \overline{\partial}^* \sigma \to \overline{\partial}^* \sigma$ in $L^2(E,e^{-2\epsilon |z|})$ as $R \to \infty$, therefore, we have $\varphi_R \sigma$ converges to $\sigma$ in the graph norm as $R \to \infty$. This finishes the proof of Claim 1. 

Next, we need to approximate (in the graph norm) any $\tau  \in $Dom$(\overline{\partial}^*)$ with compact support by elements in $C^\infty_c(E \otimes T^{1,0}\mathbb{C})$. Take any $\chi \in C^\infty_c(\mathbb{C})$ with $\int_\mathbb{C} \chi \; dxdy=1$, and set $\chi_\delta(z)=\delta^{-2} \chi(z/\delta)$. Take any $\tau  \in $Dom$(\overline{\partial}^*)$ with compact support, the convolution $\tau * \chi_\delta$ is a smooth section of $E \otimes T^{1,0}\mathbb{C}$ with compact support. Since we can fix a compact set so that all $\tau * \chi_\delta$ are supported inside the compact set, we see that $\tau * \chi_\delta$ converges to $\tau$ in $L^2(E \otimes T^{1,0}\mathbb{C},e^{-2\epsilon |z|})$ as $\delta \to 0$. It is easy to check that 
\begin{equation*} 
\overline{\partial}^* (\tau * \chi_\delta)=(\overline{\partial}^* \tau)*\chi_\delta + \epsilon \frac{\overline{z}}{|z|}(\tau * \chi_\delta) -(\epsilon \frac{\overline{z}}{|z|} \tau) * \chi_\delta. 
\end{equation*}
Therefore, 
\begin{equation*} 
\overline{\partial}^* (\tau * \chi_\delta)-(\overline{\partial}^* \tau)*\chi_\delta=  \epsilon \frac{\overline{z}}{|z|}(\tau * \chi_\delta) -(\epsilon \frac{\overline{z}}{|z|} \tau) * \chi_\delta, 
\end{equation*}
and the right hand side converges to $0$ in $L^2(E \otimes T^{1,0}\mathbb{C},e^{-2\epsilon |z|})$ as $\delta \to 0$. Hence, we conclude that $\tau * \chi_\delta$ converges to $\tau$ in the graph norm as $\delta \to 0$. Hence the proof of Lemma 6.1 is completed. 
\end{proof}

\bibliography{references}
\bibliographystyle{amsplain}

\end{document}